\documentclass[preprint,12pt,authoryear]{elsarticle}
\usepackage{graphicx} 

\usepackage{amsmath,amssymb,amsfonts,colonequals,nicefrac}
\usepackage{algorithmic}
\usepackage{graphicx}
\usepackage[rightcaption]{sidecap}
\graphicspath{{images/}} 
\usepackage{textcomp}
\usepackage{xcolor}
\usepackage{bm}
\usepackage{caption}
\usepackage{multirow}

\usepackage{siunitx}
\usepackage{url}

\usepackage{tikz}   
\usepackage{pgfplots}  
\pgfplotsset{width=7cm,compat=newest}
\usetikzlibrary{external}
\usepackage{booktabs} 
\usepackage{pgfplotstable}
\pgfplotsset{
    layers/my layer set/.define layer set={
        background,
        main,
        foreground
    }{
    },
    set layers=my layer set,
}

\newcommand{\todo}[2][]{}

\newcommand{\eg}{e.g.,\ }
\newcommand{\ie}{i.e.,\ }

\newcommand{\dx}{\mathrm{d}x}

\newcommand{\bmu}{{\bm\mu}}

\newcommand{\Mrm}{{\mathrm M}}
\newcommand{\drm}{{\mathrm d}}
\newcommand{\erm}{{\mathrm e}}
\newcommand{\srm}{{\mathrm s}}
\newcommand{\trm}{{\mathrm t}}
\newcommand{\urm}{{\mathrm u}}
\newcommand{\vrm}{{\mathrm v}}
\newcommand{\wrm}{{\mathrm w}}

\newcommand{\Rbb}{{\mathbb R}}

\usepackage{setspace}
\begin{document}
\begin{frontmatter}
\title{Optimal Experimental Design for Large-Scale Inverse Problems via Multi-PDE-constrained Optimization }
\author[1]{Andrea Petrocchi\corref{cor1}}
\ead{andrea.petrocchi@uni-konstanz.de}

\author[2]{Matthias K. Scharrer}
\ead{Matthias.Scharrer@v2c2.at}

\author[2]{Franz Pichler}
\ead{franz.pichler@v2c2.at}

\author[1]{Stefan Volkwein}
\ead{Stefan.Volkwein@uni-konstanz.de}

\cortext[cor1]{Corresponding author}

\affiliation[1]{
    organization={University of Konstanz},
    addressline={Universitaetsstrasse 10},
    postcode={78464},
    city={Konstanz},
    country={Germany},
}

\affiliation[2]{
    organization={Virtual Vehicle Research GmbH},
    addressline={Inffeldgasse 21a},
    postcode={8010},
    city={Graz},
    country={Austria},
}
\begin{abstract}
    Accurate parameter dependent electro-chemical numerical models for lithium-ion batteries are essential in industrial application.
    The exact parameters of each battery cell are unknown and a process of estimation is necessary to infer them. The parameter estimation generates an accurate model able to reproduce real cell data. The field of optimal input/experimental design deals with creating the experimental settings facilitating the estimation problem.
    Here we apply two different input design algorithms that aim at maximizing the observability of the true, unknown parameters:
    in the first algorithm, we design the applied current and the starting voltage. This lets the algorithm collect information on different states of charge, but requires long experimental times (\SI{60000}{\second}). In the second algorithm, we generate a continuous current, composed of concatenated optimal intervals. In this case, the experimental time is shorter (\SI{7000}{\second}) and numerical experiments with virtual data give an even better accuracy results, but experiments with real battery data reveal that the accuracy could decrease hundredfold.
    As the design algorithms are built independent of the model, the same results and motivation are applicable to more complex battery cell models and, moreover, to other applications.
\end{abstract}

\end{frontmatter}
\section{Introduction}
Mobile appliances and vehicles containing lithium-ion batteries have become ubiquitous over the recent decade. It is thus paramount to model and simulate lithium-ion cells efficiently and accurately. High accuracy in simulation requires the internal states of a cell to be known. These internal states may include a continuum of physical quantities, \eg potentials, concentrations and temperature, \ie a virtually infinite number of individual values. Very few quantities may actually be measured non-intrusively, \eg cell potential and temperature on the outside.

Typical high-fidelity simulation models are based on the equations described by Doyle-Fuller-Newman (DFN) model in \cite{Newman93,Newman04}. The DFN model is a homogenized system of PDEs, that typically comprise of more than 50 parameters, a vast majority of which cannot be determined experimentally.

This model may be simplified further to a Single Particle Model (SPM) by neglecting spatially resolved effects, \eg concentration differences among particles of the same electrode. Although it is known to result in less accuracy in general, it is still possible to achieve very high model accuracy by taking the SPM's limitations into account, \eg currents below \SI{1}{C}, \ie the equivalent of a full discharge in more than 1 hour, as described in \cite{Santhanagopalan06,Ning04}.

The parameter estimation for lithium-ion battery models is usually extremely time consuming. One needs very long input currents and optimization time; for example, in \cite{RSPWD19}, in order to approximate 44 battery parameters an input current longer than \SI{200000}{\second} (approximately 55 hours) is used, and even the optimization method lasts for several hours. The estimation process is based on a fitting argument: chosen a certain input current, we observe the output voltage of the battery cell. The parameter estimation is then an optimization problem based on minimizing the discrepancy between the mathematical model and the real experimental response of the battery.

Treating non-convex optimization problems, it is possible to converge to a local minimum, possibly far from the true hidden parameter. For many optimization problems this can depend on the user-defined input, whose choice change drastically the results of the parameter estimation.

In the parameter estimation framework the input current is arbitrarily chosen and, as observed in many works (\eg \cite{LWFVZ16}), the current impacts strongly on the identification of some of the parameters. For example, in \cite{LRMMJF16} an optimized battery cycling is linked to an improved parameter identifiability and therefore to robust battery health and experiment safety conditions.

The field of optimal design of experiment can help define, modify and improve the battery modeling framework. In the last decades, the literature concerning design of experiment has been extensive; we recommend some work about classic design of experiment in \cite{GP77,Paz86,Raf86,Puk06,ADT07,PP13}. Different techniques of design of experiment have been used for building new battery models, or manufacturing optimized layouts in order to improve conditions such as ageing, thermal design and charging properties; see, \eg \cite{FM08,RM22}.

The focus of this work focuses on the design of optimized inputs in the context of parameter estimation; see, \eg \cite{Meh74} for a survey on optimal input design for parameter estimation in dynamic systems. The optimal input design is strictly related to the optimal design of experiment, and the main purpose of such techniques in the framework of estimation is to improve the observability of the system.

It is well known that the many parameters in electrochemical lithium-ion battery models have different magnitude of observability -- and hence of identifiability. Some parameters are for example more sensitive to extreme \emph{states of charge}, and some are more sensitive to fast changes of applied current -- the so-called \emph{high frequency dynamics}. Standard methods of experimental design analyze the sensitivities of the parameters, namely measuring how the single parameters impact the observable quantities. In \cite{SBIG10} the parameters of a modified single-particle model are divided using a grouping algorithm and it is suggested that using experimental design methods it would be sensible to generate targeted data for the different groups.

In \cite{PKGKM18} an optimized excitation input for a parameter estimation problem is generated from a collection of predefined inputs. On the other hand, works that do not define a-priori a collection of inputs include \cite{PCVR18} and \cite{PXRS20}, where the input design of experiment for parameter identification problem is analyzed for both an isothermal and not single particle models. In particular, in \cite{PXRS20} global sensitivities are analyzed in a more rigorous approach, and a more precise parameter estimation is achieved.

In this work we start with an adaptive algorithm that, starting from some initial guesses for the parameter $\bmu^0$ and sequentially generate new input arrays $\urm^{[n]}$ at each iteration $n>0$. Each of these inputs is designed in order to maximize the observability of the current parameter, and with it we generate data and optimize the parameter $\bmu^n$; the inputs are chosen in order to improve or facilitate the parameter optimization and the convergence to the unknown underlying parameter $\bmu^*$, by collecting different pieces of information. Hence, the result of the input design algorithm is not intended to be a single ``optimal'' input, but rather a collection of inputs that collects the most information. In a later variation of the input design algorithm, we also try to optimize the experiment resources, such as minimizing the length of the excitation input. In this framework we defined optimized current sub-intervals, with ideas similar to the settings in \cite{PCVR18}.

To quantify the information of a parameter $\bmu$ carried by the current experiment, an information matrix is used, defined by using the so-called sensitivities, which measures how the output is sensitive to change in the model parameters. This information matrix resembles the Fisher information matrix, a common used statistical object and known lower bound to the covariance matrix generated by the parameter estimation. Another interpretation of the Fisher information matrix is mentioned, among others, in \cite{Uci05}, where it is seen as an approximation of the Hessian of the least-squares problem given by the parameter estimation problem. Hence, improving the well-conditioning of the approximated information matrix can also help in making the cost function more convex around the unknown parameter -- especially useful given the general non-convexity of the estimation problem.

In the standard battery cell parameter estimation applied currents are many hours long, in order to observe the battery as long as possible as collect enough information to infer the tens of parameters. Our goal is to create a collection of input current profiles much shorter than these, lasting just a few minutes, to infer 9 of the most observable parameters. One reason to treat the result of the design as a collection of inputs rather than a single ``optimal'' input is that in such a complex setting, we cannot expect to maximize the observability of all the hidden parameters with just one such short inputs, and thus we hope to capture the different information with different settings -- or inputs.

\section{Methods}
This section has the following structure: we begin in Section~\ref{Section:2.1} with the description of the numerical model. Then, we proceed to discuss the input design algorithm in Section~\ref{Section:2.2}, followed by the discretization of the design optimization and parameter estimation problems in Sections \ref{Section:2.3} and \ref{Section:2.4}.
In Section~\ref{Section:2.5} we propose a modification of the design in order to improve a weakness of the previous framework and cut the experimental time. We consider the analysis of real battery data in Section~\ref{Section:2.6} in order to study the quality of the generated optimal design.

\subsection{The numerical model}
\label{Section:2.1}
The numerical model used in this work is based on the assumption that an electrode may be fully represented by one spherical symmetric particle~-- hence it is referred to as SPM. In each particle per electrode $i\in \{A,C\}$, with $A$ and $C$ marking cathode and anode, respectively, the scaled lithium concentration $\xi_i \in [0, 1]$ 
along the respective radius ($r \in [0,R_i]$) for any time $t\in[0,t_\mathsf f]$ ($t_\mathsf f>0$) is governed by the system of equations
\begin{align}
{\dot{\xi}_i}{(t,r)}
- \frac{1}{r^2}\partial_r \big( r^2 {D_i} \partial_r {\xi_i}{(t,r)}\big)&= 0,\quad
 {(t,r) \in (0,t_\mathsf f)\times}(0,R_i),\label{eq:concentrations}
\end{align}
where $D_i>0$ is each particle's scaled diffusion constant. The spheres' symmetry is exploited by setting the natural center boundary condition at $r_i=0$ to $\partial_r {\xi_i}(t,0)=0$.
The boundary condition at the particles' outer border $R_i$ is modeled by
\begin{align}
-{D_i} \partial_r {\xi_i}(t,R_i)&=
\dfrac{1}{{\rho_i }{c_{m,i}}}{j_i}(\xi_i(t,R_i),u_i(t)),\quad{t \in (0,t_\mathsf f)},\label{eq:bc}
\end{align}
where $\rho_i c_{m,i}$ is the scaling constant of ${\xi_i}{(t,r)}$. The connection between mass transport and the electric domain is modeled in \eqref{eq:bc} by the Butler-Volmer expression
\begin{align}
{j_i}({\xi_i}{(t,R_i),u_i(t)})&={j^0_{i}}({\xi_i}{(t,R_i)}) \sinh\!\left(\tfrac{F}{R\vartheta} (u_i(t)-{O_i}({\xi_i}{(t,R_i)})\right),\label{eq:bv}
\end{align}
where $R$ is the universal gas constant and $\vartheta$ is the temperature. As the SPM formulation is expressed mainly by the open-circuit potentials $O_{i}({\xi_i}{(t,R_i)})$ of the two electrodes, these are represented by a Redlich-Kister sum as derived in \cite{Karthikeyan08}
\begin{align}
\begin{aligned}
&{O_i}({\xi_i(t,R_i)})= {U_{0,i}}+\tfrac{R\vartheta}{F} \ln\dfrac{1-{\xi_i(t,R_i)}}{{\xi_i(t,R_i)}}\\
&\qquad + \tfrac{R\vartheta}{F} \sum\limits_{k=0}^{n_i} {A_{i,k}} \bigg((2{\xi_i(t,R_i)}-1)^ {k+1}  - \dfrac{2{\xi_i(t,R_i)}k (1-{\xi_i(t,R_i)})}{(2{\xi_i(t,R_i)}-1)^{k-1}}\bigg),
\end{aligned}
\label{eq:redlich-kister}
\end{align}
where $U_{0,i}$ is the open-circuit potential offset for each electrode. Throughout the rest of the paper, we fix $U_{0,A}=0$ and use $U_{0}$ as a simplified substitute for $U_{0,C}$. The basic exchange current is modeled by 
\begin{align}
\begin{aligned}
{j^0_i}({\xi_i(t,R_i)})&= \exp\bigg(\log {k_{i}} + \tfrac{F}{R\vartheta} \big(({\xi_i(t,R_i)}-.5) O_i({\xi_i(t,R_i)})\\&\hspace{16mm}{}- \int_0^{{\xi_i(t,R_i)}} {O_i}(x)\,\dx \big) \bigg),\quad{t \in (0,t_\mathsf f)},
\end{aligned}\label{eq:models:spm:final_exchange_current}
\end{align}
where $k_i$ is the proportional reaction rate factor. For implementation, $\tilde{k}_i\colonequals\log k_i$ is used. The electric boundary condition is applied by the current $i_{\mathsf{cell}}(t)$
\begin{align}
\dfrac{3{m_i}}{{\rho_i} {R_i}}F{j_i}(\xi_i(t,R_i),u_i(t))&=-{i_{\mathsf{cell}}(t)},\label{eq:current}
\end{align}
where $F$ is Faraday's constant. The resulting cell voltage $v_{\mathsf{cell}}$ is then defined as 
\begin{align}
{u_C}(t)+{u_A}(t)+{i}_{\mathsf{cell}}(t){R_I}&={v_{\mathsf{cell}}(t)},\label{eq:voltage}
\end{align}
where $R_I$ is the total inner resistance.
\color{black}

Our work concerns the estimation of nine parameters of this model, selected on the basis of their high impact on the model behavior~-- since, as mentioned for the cases of \cite{SBIG10,PKGKM18}, different parameters are more or less identifiable.

\begin{table}[!htb]\small
    \centering
    \begin{tabular}{ cccc }
    \toprule
        Parameter & Meaning & Lower bound & Upper bound \\
        \midrule
        $D_C$ & Cathode diffusion coefficient & 0.0001 & 0.01 \\
        $D_A$ & Anode diffusion coefficient & 0.0001 & 0.01 \\
        $\xi_A$ & Anode initial state of charge & 0.007 & 0.2 \\
        $R_I$ & Inner resistance & 0.003 & 0.07 \\
        $m_C$ & Cathode active mass & 0.01 & 0.052 \\
        $m_A$ & Anode active mass & 0.006 & 0.034 \\
        $\tilde{k}_C$ & Cathode reaction rate & -20 & -1 \\
        $\tilde{k}_A$ & Anode reaction rate & -23 & 10 \\
        $U_0$ & Cathode 0th RK parameter & 3 & 4 \\
        \bottomrule
    \end{tabular}
    \caption{Unscaled parameter names and bounds.}
    \label{2:table_unscaled_parameters}
\end{table}

The collection of considered parameters is depicted in Table~\ref{2:table_unscaled_parameters}, together with lower and upper bounds. In particular, the diffusion coefficients $D_C, D_A$ appear in \eqref{eq:concentrations} as well as the initial state of charge of the anode $\xi_A$~-- since $\xi_C$ can be recovered algebraically from a given cell voltage $v\vert_{t=0}$. The total inner resistance $R_I$ appears in the voltage equation \eqref{eq:voltage}, while the masses $m_C, m_A$ are involved in the current equation \eqref{eq:current}. The reaction rate logarithms $\tilde{k}_C, \tilde{k}_A$ are relevant in the exchange current formula \eqref{eq:models:spm:final_exchange_current} and, finally, the first Redlich-Kister term $U_0$ is considered, while all the other model parameters are set to experimental acceptable standards.

As these values have different cardinalities, we apply the following scalings:
\begin{align*}
	& \mu_1 := \log_{10} \left( \frac{D_C}{D_C^L} \right), &
	& \mu_2 := \log_{10} \left( \frac{D_A}{D_A^L} \right), &
	& \mu_3 := \frac{\xi_A}{(\xi_A^L+\xi_A^U)/2}, \\
	& \mu_4 := \frac{R_I}{(R_I^L+R_I^U)/2}, &
	& \mu_5 := \frac{m_C}{(m_C^L+m_C^U)/2}, &
	& \mu_6 := \frac{m_A}{(m_A^L+m_A^U)/2}, \\
	& \mu_7 := \tilde{k}_C - (\tilde{k}_C^L+\tilde{k}_A^U)/2, &
	& \mu_8 := \tilde{k}_A - \tilde{k}_C, &
	& \mu_9 := \frac{U_0}{(U_0^L+U_0^U)/2}.
\end{align*}

Finally, we report the scaled parameter bounds in Table~\ref{2:table_scaled_parameters}.
\begin{table}[!htb]\small
    \centering
    \begin{tabular}{ cc }
    \toprule
        $\bmu^L$ & $\bmu^U$ \\
        \midrule
        0	&	2	\\
        0	&	2	\\
        0.0676328502415459	&	1.93236714975845	\\
        0.0821917808219178	&	1.91780821917808	\\
        0.32258064516129	&	1.67741935483871	\\
        0.3	&	1.7	\\
        -9.5	&	9.5	\\
        -3	&	11	\\
        0.857142857142857	&	1.14285714285714	\\
        \bottomrule
    \end{tabular}
    \caption{Bounds for the scaled parameters}
    \label{2:table_scaled_parameters}
\end{table}

\subsection{The input design algorithm}
\label{Section:2.2}
From now on, we use the battery model as a black box system that takes a model parameter vector $\bmu \in \Rbb^d$ and a current density function $i:[0,t_\mathsf f  ] \to \Rbb$ input variables, whereas the voltage $v_{\bmu,i}: [0,t_\mathsf f] \to \Rbb$ is the output variable.

We use an adaptive algorithm that starts from an arbitrary parameter $\bmu^0$ and, for $n\ge1$ (called \emph{design iteration}), the algorithm iterates on three phases:
\begin{itemize}
    \item In the first phase we generate a new current function $i^{[n]}$.
    \item In the second phase we read and store some data $w^{[n]}$ derived by using the current $i^{[n]}$.
    \item The last phase is the parameter optimization/estimation using previous data $\{w^{[1]},\dots,w^{[n]}\}$, with which we obtain the parameter $\bmu^{n}$, approximation of the true hidden parameter $\bmu^*$.
\end{itemize}

There are different ways of ending this algorithm. One possible way is to stop the algorithm as soon as the last generated current function is sufficiently close to one of the already computed currents $i^{[k]}$, $k=1,\ldots,n-1$, in our collection. Therefore, we use the stopping criterion
\begin{equation}
	\label{4.2:OID_stopping_criterion_tol}
    \min_{i \in \{ i^{[1]},\ldots,i^{[n-1]}\} } \| i^{[n]} - i \|_{L^2(0,t_\mathsf f)} < \varepsilon
\end{equation}
with a properly chosen tolerance $\varepsilon>0$. However, a known limitation is that we have no guarantee of convergence. This fact can potentially generate an infinite number of current inputs. In the cases, where we have used this criterion, we have also set a maximum number of iterations as a safeguard.
An alternative termination is to fix \emph{a priori} a determinate number of inputs and forcing the algorithm to generate inputs sufficiently far from each other, by penalizing the closeness between them. We will later describe in detail this forced generation of inputs.

Next, we introduce a so-called symmetric information matrix in order to quantify the uncertainty produced by a certain input $i$. As in \cite{ADT07}, we define the information matrix $\Mrm \in \Rbb^{d\times d}$ by using the sensitivity variables, described in details in the next subsection. Then, the uncertainty given by an information matrix is defined as the inverse by applying the D-optimal design, or D-criterion (see, \eg \cite{GP77,ADT07}), \ie
\begin{equation*}
    \phi(\Mrm) \colonequals - \log \left( \det \left( \Mrm \right) \right),
\end{equation*}
which is equivalent to minimizing the generalized variance of the parameters -- as discussed, \eg in \cite{ADT07}.
We can also define the uncertainty function by using the eigenvalues of $\Mrm$ -- especially convenient when the dimension $d$ is high and evaluating the determinant gets computationally expensive. Indeed, if $\Lambda := \{ \lambda_1,\dots,\lambda_d\}$ is the spectrum of $\Mrm$, then
\begin{equation*}
	\phi(\Mrm) \colonequals \log \left( \prod_{\lambda_j \in \Lambda} \frac{1}{\lambda_j} \right).
\end{equation*}
Observe that we can also interpret this criterion as minimization of the volume of the confidence ellipsoid for the parameters.

\subsection{Approximation of the uncertainty measure}
\label{Section:2.3}
We fix the time step to $\delta=0.1$ seconds, so that we have $K=600$ time steps in an interval of \SI{60}{\second}. 
Time steps $t_k:=k\delta$ for $k=0,\dots,K$ are discretized uniformly, and all current functions will be discretized as step function evaluated on such time steps, so that the current information information is compressed into the \emph{input array} $\urm$ -- or, simply, \emph{input}. Given $n_\urm \in \mathbb{N}^+$ and $\tau_\urm \colonequals t_\mathsf f / n_\urm$, the array
\begin{equation}
	\urm \colonequals [ u_1, \dots, u_{n_\urm} , v_0 ] \in \Rbb^{n_\urm+1}
\end{equation}
corresponds to the current function
\begin{equation*}
    i_\urm (t) \colonequals 
    \left\{
    \begin{aligned}
        &u_1&&\text{for }0\le t<\tau_\urm,\\
        &u_2&&\text{for }\tau_\urm \le t<2\tau_\urm,\\
        &&&\vdots\\
        &u_{n_\urm}&&\text{for }(n_\urm-1) \tau_\urm \le t \le t_\mathsf f,
    \end{aligned}
    \right.
\end{equation*}
and the last component, $v_0$, is the initial voltage of the battery, which is required by the numerical model. With $v_{\bmu,\urm} : [0,t_\mathsf f] \to \mathbb{R}$ we indicate the output voltage given the parameter $\bmu$ and input array $\urm$ (hence given current $i_\urm$ and initial voltage $v_0$). We seldom use its discretization $\vrm_{\bmu,\urm}$ evaluated on the time steps $\{t_k\}_{k=0}^K$.

The sensitivity variables (or, simply, \emph{sensitivities}) are usually evaluated by solving the so-called \emph{sensitivity equations}. However, in a black-box model they can only be approximated by using finite differences. For $1 \le j \le d$, the $j$-th discretized sensitivity $\srm^j=\srm^j_{\bmu,\urm} \in \Rbb^{K+1}$ is defined as
\begin{equation*}
    \srm^j_{\bmu,\urm} := \frac{\vrm_{\bmu + \nu \erm_j,\urm} - \vrm_{\bmu,\urm}}{\nu} \approx \frac{\drm}{\drm \mu_j} \vrm_{\bmu,\urm},
\end{equation*}
where $\erm_j$ is the $j$-th Euclidean basis element, and $\nu:=10^{-3}$. Let us mention that the approximation of the sensitivities can be improved by using other (and more accurate) methods to evaluate the partial derivatives, such as automatic differentiation (see, \eg \cite{Ra81}).

Then, the information matrix $\Mrm=\Mrm_\bmu^\urm$, evaluated at some parameter $\bmu$ and input array $\urm$, is defined as
\begin{equation*}
    \left( \Mrm \right)_{j,l} \colonequals \sum_{k=0}^{K-1} \frac{ \srm^j_k \srm^l_k + \srm^j_{k+1} \srm^l_{k+1} }{2} ( t_{k+1} - t_k ), \qquad 1\le j,l \le d,
\end{equation*}
namely, as the approximation of the $L^2$ inner product $\langle s^j, s^l \rangle_{L^2(0,t_\mathsf f)}$ via the trapezoidal rule, where the sensitivities are evaluated at the parameter $\bmu$ and the current $i_\urm$.

As we want to maximize the total information coming from a collection of inputs, we modify the information matrix of each design iteration. In particular, we simply sum the information matrices, thus at design iteration $n$ we define the \emph{optimal design function}
\begin{equation*}
	\Phi(\urm) \colonequals - \log_{10} \left( \det \bigg( \sum_{m=1}^{n-1} \Mrm_{\bmu^{n-1}}^{[m]} + \Mrm_{\bmu^{n-1}}^\urm \bigg) \right) + \gamma\,{\| \urm \|}^2_2,
\end{equation*}
where $\Mrm^{[m]}$, $m=1,\dots,n-1$, denotes the information matrix corresponding to input array $\urm^{[m]}$, $\gamma\colonequals 10^{-4}$ accompanies a regularization term and $\|\cdot\|_p$ indicates the $\ell^p$-norm for $1\le p \le \infty$.

Hence, to find a new input array we solve the optimization problem
\begin{equation}
    \label{3.1:design_optimization_problem}
	\min_{\urm \in \Rbb^{n_\urm+1}} \Phi(\urm)\quad \text{subject to}\quad \urm^{\mathsf{min}} \le \urm \le \urm^{\mathsf{max}},
\end{equation}
where $\urm^{\mathsf{min}}$ and $\urm^{\mathsf{max}}$ are the bounds in $\Rbb^{n_\urm+1}$ chosen as
\begin{align*}
    & -8.8 \equalscolon u^{\mathsf{min}}_j \le u_j \le u^{\mathsf{max}}_j \colonequals 8.8,\quad j=1,\dots,n_\urm,\\
    & 3.3 \equalscolon u^{\mathsf{min}}_{n_\urm+1} \le v_0 \le u^{\mathsf{max}}_{n_\urm+1} \colonequals 4.1.
\end{align*}
in our experiments, where the unit of current is \unit{\ampere} and the unit of the initial voltage is \unit{\volt}.

The optimization problem \eqref{3.1:design_optimization_problem} is solved numerically using the routine {\texttt{fmin\_l\_bfgs\_b}} from the \texttt{scipy} Python library; see, \eg \cite{scipy}.

As mentioned earlier, when we want the algorithm to generate a fixed number of inputs, we need to make them distinct, since there is no information gain when the input generated is equal to one we already have. We therefore define an alternative (nonsmooth) optimization function that enforces the new inputs to be sufficiently far from each other
\begin{equation}
    \label{3.1:uncertainty_with_penalization}
    \hat\Phi(\urm) \colonequals \Phi(\urm) + {\textstyle \sum\limits_{\substack{\urm^{[n]} \text{ previous}\\\text{controls}}}} \frac{1}{1 + 100 \, {\|\urm - \urm^{[n]} \|}_\infty}.
\end{equation}
The penalizing function, as well as the scalar 100 have been chosen in an intuitive way in order to guarantee that we stop our algorithm provided our current control is sufficiently far from the already computed controls, but also not to modify the initial optimization problem too strongly.

\subsection{Parameter estimation}
\label{Section:2.4}
Assuming now that the input vector $\urm$ (hence the current $i$) is fixed, we want to find the parameter vector $\bmu$ minimizing the discrepancy between our parametrized model and data coming from the real application. We assume, for the moment, that the data is given by the voltage corresponding to some unknown parameter $\bmu^*$, evaluated on the time array $\trm$, \ie
\begin{equation*}
    w_\urm = v_{\bmu^*,\urm} : [0,t_\mathsf f] \to \mathbb{R},
\end{equation*}
and its discretization is indicated as $\wrm_\urm \colonequals \vrm_{\bmu^*,\urm} \in \Rbb^{K+1}$. This means that we consider noiseless virtual data.

For the numerical optimization we use a least squares method (from \texttt{scipy.optimize}, see \cite{BCL99}) to fit $\vrm_{\bmu,\urm}$ to $\wrm_\urm$. In particular, we minimize the relative error $\erm_{\urm}(\bmu) \colonequals (\vrm_{\bmu,\urm}-\wrm_\urm)/\wrm_\urm$ (where the division is intended component-wise) by in particular solving
\begin{equation}
    \label{3.2:least_squares_problem}
    \min_{\bmu\in \mathcal{P}_{\mathsf{ad}}} J(\bmu) \colonequals \frac{1}{2} {\left\| \frac{\vrm_{\bmu,\urm} - \wrm_\urm}{\wrm_\urm} \right\|}^2_2,
\end{equation}
where $\mathcal{P}_{\mathsf{ad}} \subset \mathbb{R}^d$ is the closed, convex set of admissible parameters -- in our case the box-constrained set defined by $\bmu^L \le \bmu \le \bmu^U$ in Table~\ref{2:table_scaled_parameters}.

In the framework of multiple data, we can simply merge information by stacking the data vectors together. In particular, given $\urm^{[1]}, \ldots, \urm^{[n]}$ input arrays generated by the algorithm, and $\wrm^{[1]},\ldots,\wrm^{[n]}$ associated data vectors (with $\wrm^{[n]}\colonequals\wrm_{\urm^{[n]}}=\vrm_{\bmu^*,\urm^{[n]}}$), then we want to fit $\vrm_{\bmu,\urm^{[1]}}$ to $\wrm^{[1]}$ and so on.
We concatenate all vectors of relative error into one vector $[\erm_{\urm^{[1]}}(\bmu),\ldots,\erm_{\urm^{[n]}}(\bmu)] \in \Rbb^{n(K+1)}$ and we minimize its $\ell^2$-norm.
The parameter generated at OID iteration $n$ is then the optimizer $\bmu^n$.
We mention that in order to optimize computing time, we evaluate voltage vectors in a parallel way for all different current inputs.

We later explore the effect of the input design on the parameter optimization problem by analyzing the Hessian matrix $H$ of the least-squares function $J$ defined in \eqref{3.2:least_squares_problem}, evaluated in the optimum $\bmu^*$. In particular, we analyze the value
\begin{equation}
    \label{3.2:beta}
    \beta \colonequals \lambda_{\mathsf{max}}(H)/\lambda_{\mathsf{min}}(H),
\end{equation}
i.e., the conditioning number of the Hessian, where at iteration $n$ this value will be indicated as $\beta_n$ -- as the Hessian matrix changes at each OID iteration $n$. When this value decreases, the smallest eigenvalue is growing because in this construction the biggest eigenvalue cannot decrease, hence indicating positivity of the Hessian and convexity of the optimization problem.
Let us comment on the fact that we cannot use the term $\beta$ in the definition of the optimal input design optimization problem, since this quantity depends on knowing the hidden parameter $\bmu^*$. Nonetheless, we hope to minimize this value indirectly.

\subsection{Another design framework}
\label{Section:2.5}
Generating ten different short inputs sounds like a good idea if we want to look at the behavior of a battery cell in all different states of charge: starting from different initial voltages serves exactly this purpose. But there might be a problem in the actual use of these different current profiles, since we cannot evaluate all inputs at the same time in a parallel way -- as we do in the numerical model. On the contrary, before and after every input we need a phase of preparation to bring the battery cell at the right initial voltage, and a long phase of rest to let the internal electro-chemical dynamics stabilize. For this reason, even if ten different inputs of \SI{60}{s} each together sum up to \SI{600}{s}, actually setting up and running the experiment can take more than \SI{60000}{s}.

Hence, we have modified the design framework and structure in order to avoid this problem. The idea is to generate a sequence of new input intervals successively, rather than to use them as entirely different current profiles. As a consequence, we may not choose the initial voltage of each interval, because this will be given by the ending voltage of the previous interval. Rather, we may only choose the initial voltage, namely the starting voltage of the first interval; this is chosen as \SI{3.9}{\volt}, equivalent to appr.\ \SI{75}{\percent} SOC.

\begin{figure}[!htb]
	\centering
	\includegraphics{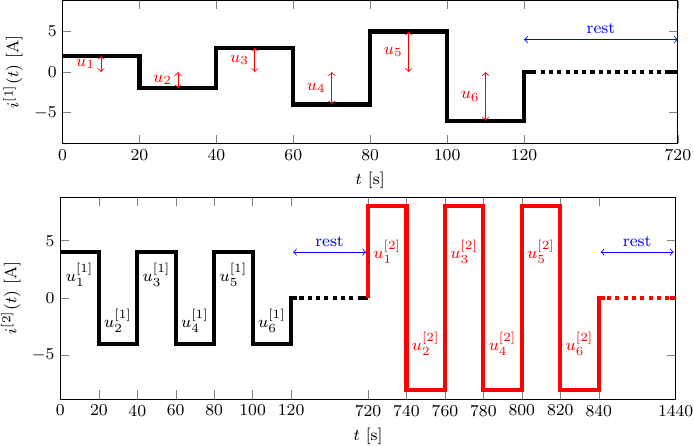}
	\caption{Visual arbitrary example of the concatenated-inputs setting. Above: initial current function, where the input variables define the 6 jumps and a successive rest phase is imposed. Below: once the first current function is fixed, the next input array $\urm^{[2]}$ is optimized by optimizing the concatenated current function $i^{[2]}$.}
    \label{4:figure_continuous_example}
\end{figure}

We also realize that we need a resting interval in order to capture the long-term dynamics. Hence, our input is translated into six ``jumps'' of a \SI{120}{\second} high-frequency time interval followed by a \SI{600}{\second} resting phase. An example of the first input with arbitrary jumps is plotted in Figure~\ref{4:figure_continuous_example}.
This time we are not arbitrarily choosing the first input $\urm^{[1]}=(u^{[1]}_1,\dots,u^{[1]}_6)$, but it is found in the first step of the OID algorithm.
In the following iterations, we find input array $\urm^{[n]}$ by optimizing the uncertainty of the full-concatenated current input, i.e., we add the new input to the old one and hence consider only one longer input. See Figure~\ref{4:figure_continuous_example} for a visual explanation.

\subsection{Analysis on real battery data}
\label{Section:2.6}
While we evaluate these optimal inputs using virtual data $w_\urm=v_{\bmu^*,\urm}$ in the algorithm, we want to check how beneficial they are in a real-world parameter estimation problem. To this end, we have run experiments on a real battery cell and collected the measurement data.
The measured voltage data obtained through cycling two Molicel INR21700-P45B commercial cells is used in this work.
For a given starting voltage, the cell is first fully charged at a constant current of \SI{4.5}{A} and then, at reaching the top voltage of \SI{4.2}{V}, switch to constant voltage mode until the magnitude of the current tapered down below \SI{0.05}{A}. After a resting period of \SI{600}{s}, the cell is discharged at \SI{4.5}{A} to the given voltage followed again by a constant voltage step until \SI{0.05}{A} is underrun.
Then, the designed current profile is applied in a ``current simulation'' step.
All tests were carried out using an Arbin LBT21084HC battery testing system and a Memmert incubator with Peltier cooling (model IPP600) for maintaining the temperature at \SI[separate-uncertainty = true]{25(2)}{\celsius} by forced air cooling.%

We want to use these real-world data to measure the effect of the optimal input design on the parameter estimation. This data also includes preparation steps at the beginning and the resting phase at the beginning, and depending on the test we will use information coming from the extra phases or not. Furthermore, the time stepping was chosen depending on the situation, varying between a minimum of \SI{0.0002}{s} (used in high-frequency dynamics) and a maximum of \SI{1.0}{s} (used in the resting phases).

\section{Results and Discussion}
\subsection{First input design: collection of inputs}
\label{Section:3.1}
\begin{figure}[!ht]
    \centering
    \includegraphics[width=\textwidth]{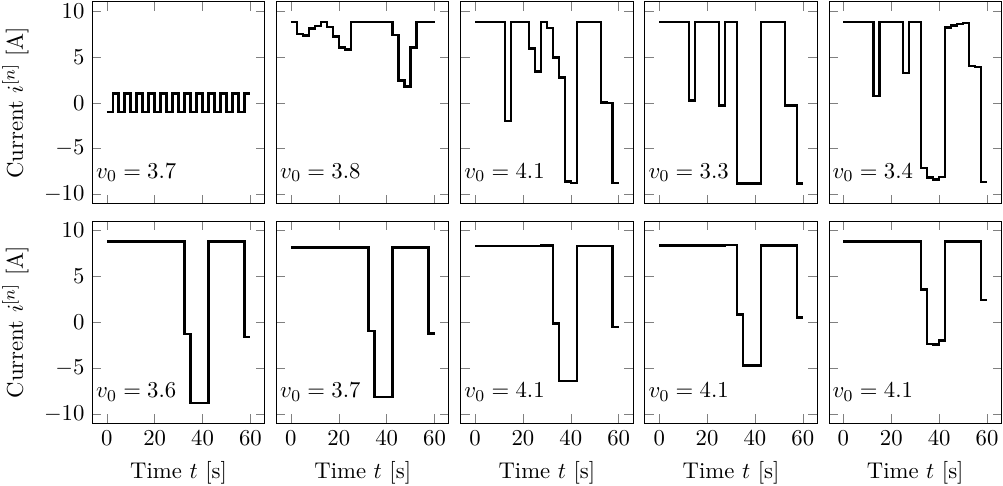}
    \caption{Input functions $i^{[n]}$, $n=1,\dots,10$ generated by the algorithm. In the lower left corner the initial voltages $v_0$ are indicated.}
    \label{3.3:figure_adaptive_inputs}
\end{figure}

We have applied the input design algorithm with piecewise constant inputs consisting of 24 jumps, therefore each continuous sub-interval is \SI{2.5}{\second} long and the input dimension is 25. We furthermore set the maximum number of inputs to 10 and arbitrarily fix the first input as
\begin{equation}
	\label{3.3:adaptive_initial_input}
	\urm^{[1]} \colonequals [-1, +1, -1, +1, \dots, -1,+1, 3.7] \in \Rbb^{25}.
\end{equation}
This corresponds to a quickly alternating input applied to a roughly mid-charge battery. The generated inputs are plotted in Figure~\ref{3.3:figure_adaptive_inputs}.

Notice how some of the inputs look very similar, but the penalization function introduced in \eqref{3.1:uncertainty_with_penalization} guarantees that we do not use the exact same input array more than once.

\begin{figure}[hp]
	\centering
	\includegraphics[scale=0.8]{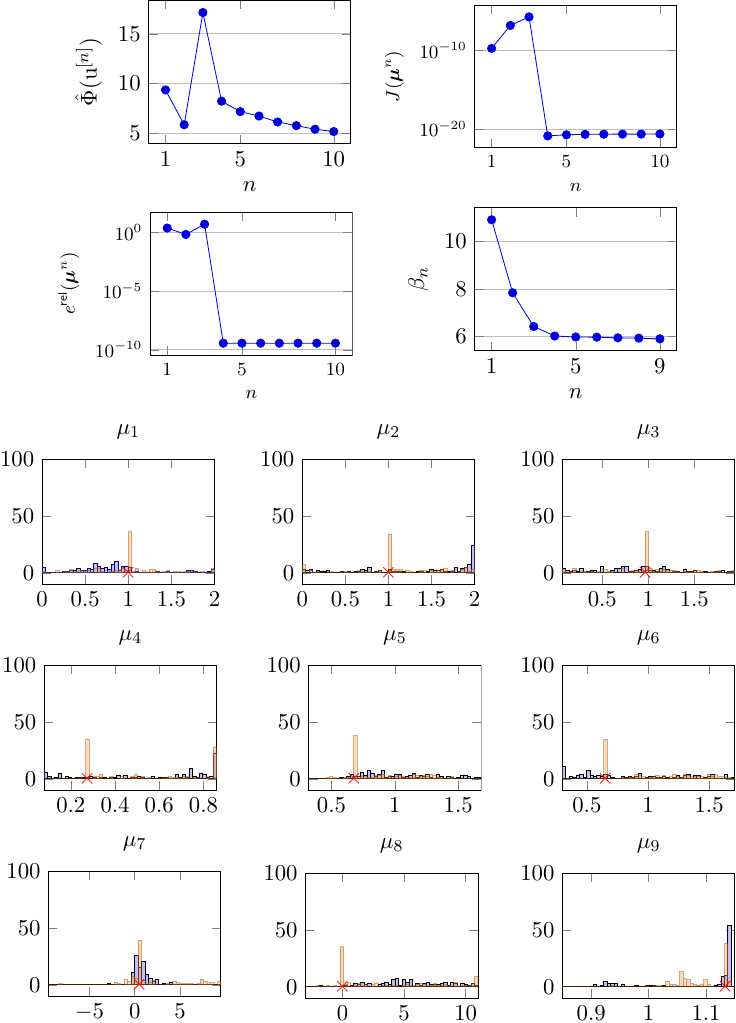}
	\caption{Top, blue: information corresponding to each design iteration $n$. First line: uncertainty function and cost evaluated in the optimal parameter. Second line: Euclidean norm of the error between the optimal parameter and the true parameter and value $\beta_n$ given in \eqref{3.2:beta}. Bottom: qualitative statistical study for the input design algorithm. Dark bars (blue) show the results of 100 parameter estimations using only the input given by $\urm^{[1]}$. Light bars (orange) show the results using the optimal collection of inputs. The red crosses mark the hidden parameters.}
	\label{3.3:figure_adaptive_analysis}
\end{figure}

In order to analyze the effect of the input design algorithm, in Figure~\ref{3.3:figure_adaptive_analysis} we plot the change in different quantities, such as the design algorithm the uncertainty function $\hat\Phi$, the parameter optimization cost $J$, the relative error $e^{\mathsf{rel}}(\bmu^n) := \|\bmu^*-\bmu^n\|_2 / \| \bmu^*\|_2$ and the sequence $\beta_n$ defined in \eqref{3.2:beta}.
In the same plot, in the lower area we plot a qualitative analysis of the input design, where we compare the results of 100 parameter optimizations with randomized initial parameter at the beginning and at the end of the algorithm, \ie comparing the information before (in dark colored bars) and after (in light colored bars).
Firstly, we can comment on the properties of the initial input $\urm^{[1]}$ defined in \eqref{3.3:adaptive_initial_input}: while it looks like the error $J$ is already low ($1.77\cdot10^{-10}$) at the initial parameter, the relative error of the parameter reads in fact 2.22.
Namely, the initial input holds little information and it does not provide data such that the parameter optimization approximates $\bmu^*$. We can observe this in the bar plot, too: the blue bars do not approximate the parameter indicated with the red cross most of the time.

During the design algorithm the uncertainty function, the parameter optimization cost and the relative error do not decrease monotonously. The OID algorithm converges to the true parameter after some preliminary steps, in particular reaching the neighborhood of $\bmu^*$ in the fourth iteration.
This could be explained by the behavior of the sequence $\beta_n$, which decreases monotonously, reaching a plateau after four iterations, when the algorithm finds the hidden parameter with a relative error of approximately $3.73 \cdot 10^{-10}$.

In the bar plots of Figure~\ref{3.3:figure_adaptive_analysis} we can observe that the input design improves the parameter estimation significantly: with a non-informative input (dark bars), the parameter estimation fails almost every time, while using the collection of selected inputs (light  bars) the convergence to the hidden parameter (the red cross) improves greatly.

Although the converge to the optimal parameter does not happen always, the optimization results cluster around the optimal parameters. Overall, we can see a degree of bias in the bar plots of Figure~\ref{3.3:figure_adaptive_analysis} -- namely, we can guess the hidden parameter in an easier way. However, this framework might find in its applicability some unexpected counter-effects, which we discuss in the next section.

\subsection{Second input design: concatenated input}
\label{Section:3.2}
\begin{figure}[!htb]
	\centering
	\includegraphics[scale=0.8]{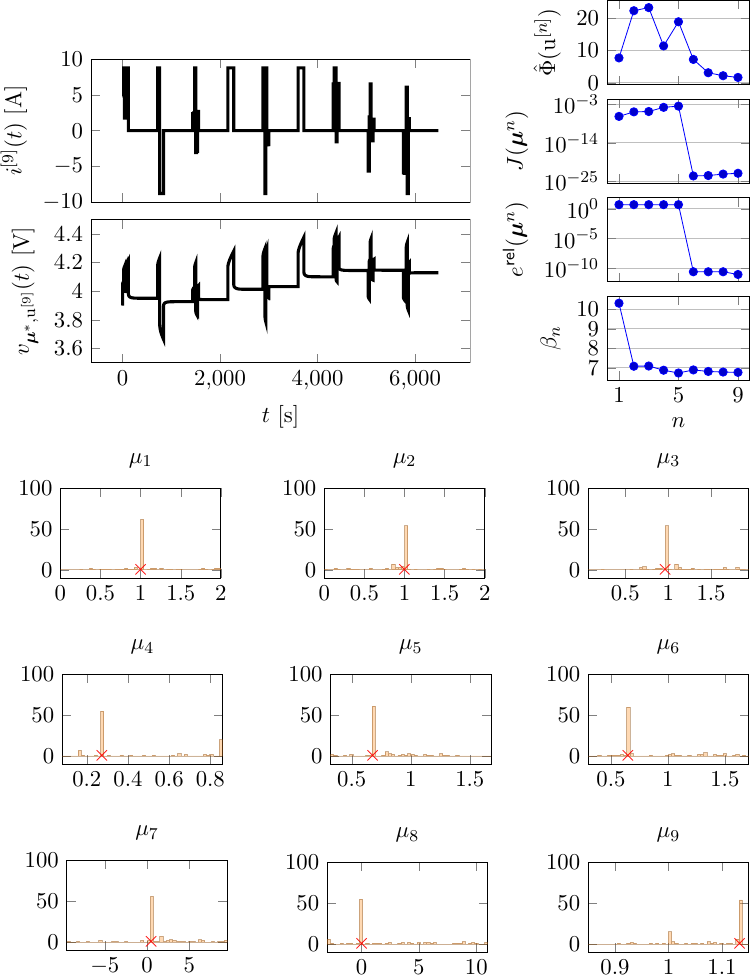}
    \caption{Top left, in black: optimal input generated by the algorithm and corresponding output voltage generated by the battery cell model. Top right, in blue: plots of different quantities at OID iteration $n$. Bottom bar plots: qualitative statistical results of the optimal parameter of 100 parameter estimations. The red crosses mark the true parameters $\bmu^*$.}
    \label{4:figure_continuous_results}
\end{figure}

In this setting we generate 9 inputs, each lasting \SI{120}{\second}: therefore, concatenating the respective jumps and rest phases together we reach an ``optimal'' input lasting \SI{6480}{\second}, i.e. 108 minutes. This optimal current $i^{[9]}$ generated by the algorithm is plotted in Figure~\ref{4:figure_continuous_results} together with the corresponding output voltage $v_{\bmu^*,\urm^{[9]}}$ produced from the battery cell model using such input.
Further information for each OID iteration are also plotted in Figure~\ref{4:figure_continuous_results}. In this case we can observe a non monotone decrease of the quantity $\beta_n$, but still a good convergence of the parameter optimization in later iterations, reaching a relative parameter error of $9.74 \cdot 10^{-12}$.

Similar to what we did previously, in the bottom part of Figure~\ref{4:figure_continuous_results} we plot the result of 100 separate executions of the parameter optimization in order to show the estimation error for randomized initial parameters.
We observe an overall better consistence and accuracy in the convergence to the true parameter; considering the fact that we the overall experiment time will be reduced at least 4 times, as shown in the next section, we regard this as a great improvement of the input design presented in the previous section.

\subsection{Real battery data}
\label{Section:3.3}
We consider four different situations, depending on two different factors: whether we use the collection of inputs derived in Section~\ref{Section:3.1} or the concatenated input derived in Section~\ref{Section:3.2}, and whether we consider the additional data provided by the preparation and resting phases, or not. These choices are summarized in Table~\ref{5:different_tests}, and some of the results given in each test are summarized in Table~\ref{5:tests_results}.

\begin{table}[!ht]
    \centering
    \begin{tabular}{ ccc }
        \toprule
        & Cut data: & Full data:\\
        & only impulses & impulses + rest phases \\
        \midrule
        Collection of inputs & Test 1 & Test 2 \\
        Concatenated input & Test 3 & Test 4 \\
        \bottomrule
    \end{tabular}
    \caption{Considered tests for the analysis of real data.}
    \label{5:different_tests}
\end{table}

\begin{table}[!htb]
    \centering
    \begin{tabular}{ ccccc }
        \toprule
        & \multicolumn{4}{c}{TESTS}\\
        & 1 & 2 & 3 & 4 \\
        \midrule
        Experiment time ($s$) & 600 & 63,362 & 6,480 & 14,199 \\
        Data points & 20,350 & 83,160 & 216,002 & 223,727 \\
        Optimization time ($s$) & 800 & 154 & 550 & 450 \\
        \bottomrule
    \end{tabular}
    \caption{First line: total input time $t_\mathsf f$; second line: number of data points; third line: average time for the parameter optimization.}
    \label{5:tests_results}
\end{table}

In each of these situations we run 100 parameter optimizations in order to infer the unknown battery cell parameters. We analyze the quality of the estimation using a bar plot, and we try to infer the true parameter by proposing an optimal parameter. Since this is unknown, we finally perform a qualitative error analysis on the proposed parameters given in each test.

In particular, we will propose three parameters $\bmu^{\mathsf{opt},2}$, $\bmu^{\mathsf{opt},3}$ and $\bmu^{\mathsf{opt},4}$ inferred from test 2, 3 and 4, respectively. For each of these parameters we plot the output voltage produced by the numerical model $v_{\bmu,\urm}(t)$, as well as the relative error given by the discrepancy with the real battery data $w(t)$, namely
\begin{equation*}
	e^{\mathsf{rel},[n]}_{\boldsymbol\mu^{\mathsf{opt},k}} (t) := (w^{[n]} (t) - v_{\bmu^{\mathsf{opt},k}, \urm^{[n]}} (t)) / w^{[n]} (t) \qquad \text{for } n=2,3,4 \ \text{and } k=2,3,4.
\end{equation*}
The results are shown in the upcoming figures.

\subsubsection*{Test 1}
\begin{figure}[!htb]
    \centering
    \includegraphics[scale=0.8]{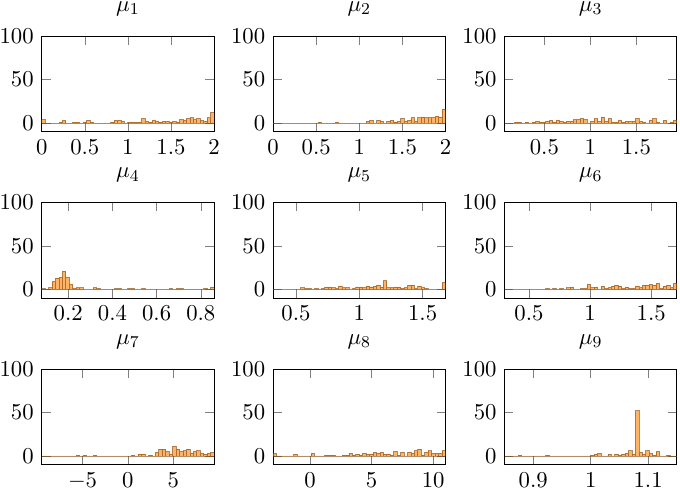}
    \caption{Results for Test 1. Bar plots of qualitative analysis of 100 parameter estimations starting from randomized initial guess.}
    \label{5:figure_test_1}
\end{figure}

By counting only the 10 one-minute intervals, we get a relatively small experiment time and number of data points. Although the optimization is longer than in other cases, we cannot identify the unknown parameters with confidence. Even if we are using the same input collections as the results shown in Figure~\ref{3.3:figure_adaptive_analysis}, the bar plots shown in Figure~\ref{5:figure_test_1} look much worse. Therefore, in this case we conclude that this input design is not successful, as we derive no hint about where the hidden parameters might be.

\subsubsection*{Test 2}
\begin{figure}[!htb]
    \centering
    \includegraphics[scale=0.8]{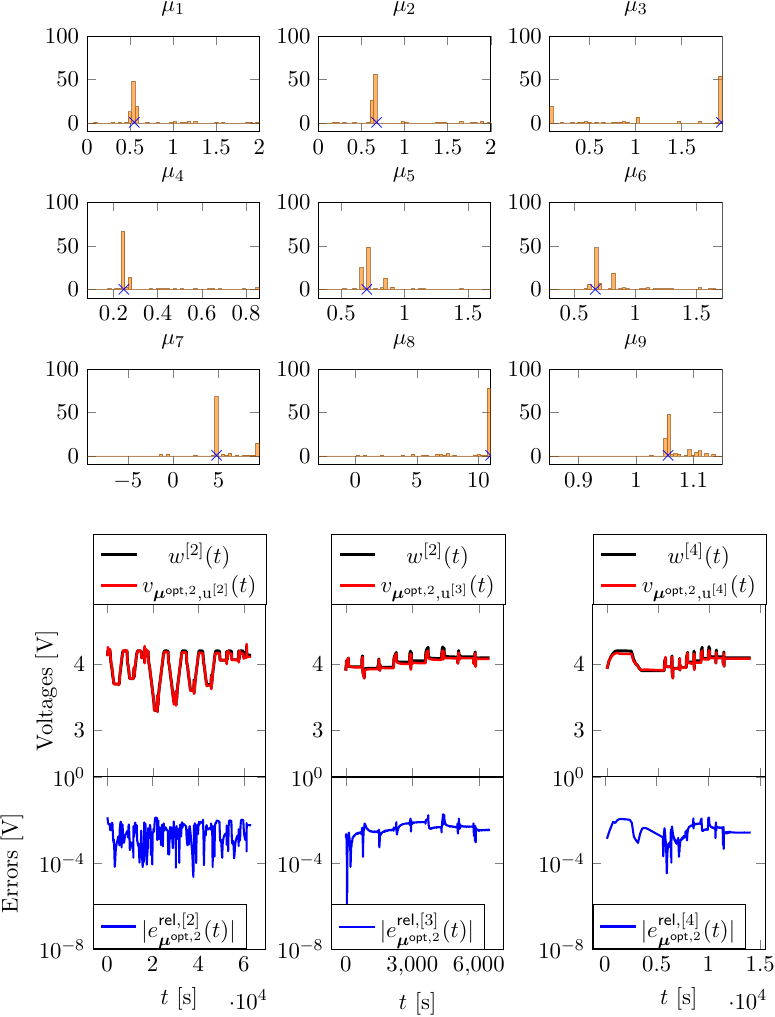}
    \caption{Results for Test 2. Top: bar plots of qualitative analysis of 100 parameter estimations starting from randomized initial guess; the blue cross finds the proposed parameter $\bmu^{\mathsf{opt},2}$; bottom: data, output and error generated by the proposed parameter.}
    \label{5:figure_test_2}
\end{figure}

Since the production of the data necessary for Test 1 produces, by default, the data used in Test 2, this test is necessarily more informative and hence the estimation of the true parameter is more successful. However, the experiment lasts for over 17 hours, with an average optimization time of about two and a half minutes. As shown in Figure~\ref{5:figure_test_2}, we can identify a candidate $\bmu^{\mathsf{opt},2}$, indicated with the blue mark.

The analysis of the model error given by the proposed parameter confirms the quality of such parameter: for the inputs $\urm^{[2]}$ and $\urm^{[3]}$ the relative error oscillates in the order of $10^{-3}$, while for the input $\urm^{[4]}$ the order of the relative error is $10^{-2}$. The numerical model tuned with such parameter is therefore an accurate simulator of the battery cell taken into consideration.

\subsubsection*{Test 3}
\begin{figure}[!htb]
    \centering
    \includegraphics[scale=0.8]{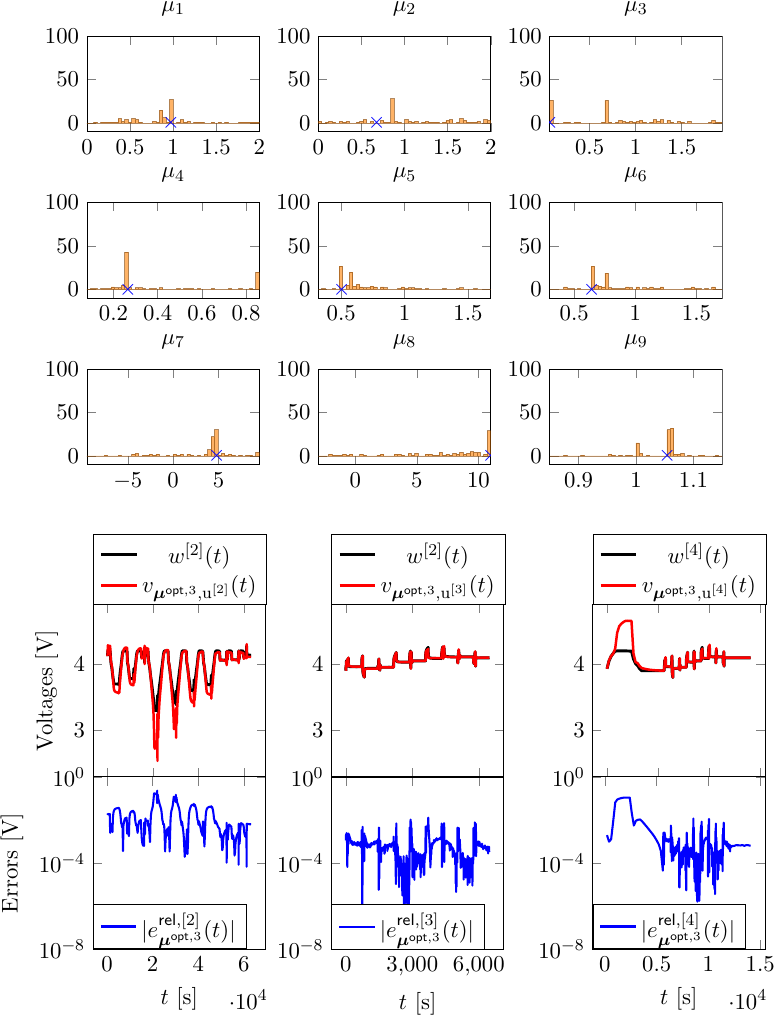}
    \caption{Results for Test 3. Top: bar plots of qualitative analysis of 100 parameter estimations starting from randomized initial guess; the blue cross finds the proposed parameter $\bmu^{\mathsf{opt},3}$; bottom: data, output and error generated by the proposed parameter.}
    \label{5:figure_test_3}
\end{figure}

By construction, the time of the generated experiment is 108 minutes, and the optimization over 216,002 points takes in average several minutes.
However, we can see that in this case the estimations gravitate around different local minima. With $\bmu^{\mathsf{opt},3}$ we have indicated the most common optimizer, indicated in Figure~\ref{5:figure_test_3} with the blue mark.

This parameter produces accurate results for the data $w^{[3]}$, but accumulates a too high error in the higher voltage changes occurring in data observations $w^{[2]}$ and $w^{[4]}$, hence not capturing very well the full range of battery cell dynamics.

\subsubsection*{Test 4}
\begin{figure}[!htb]
    \centering
    \includegraphics[scale=0.8]{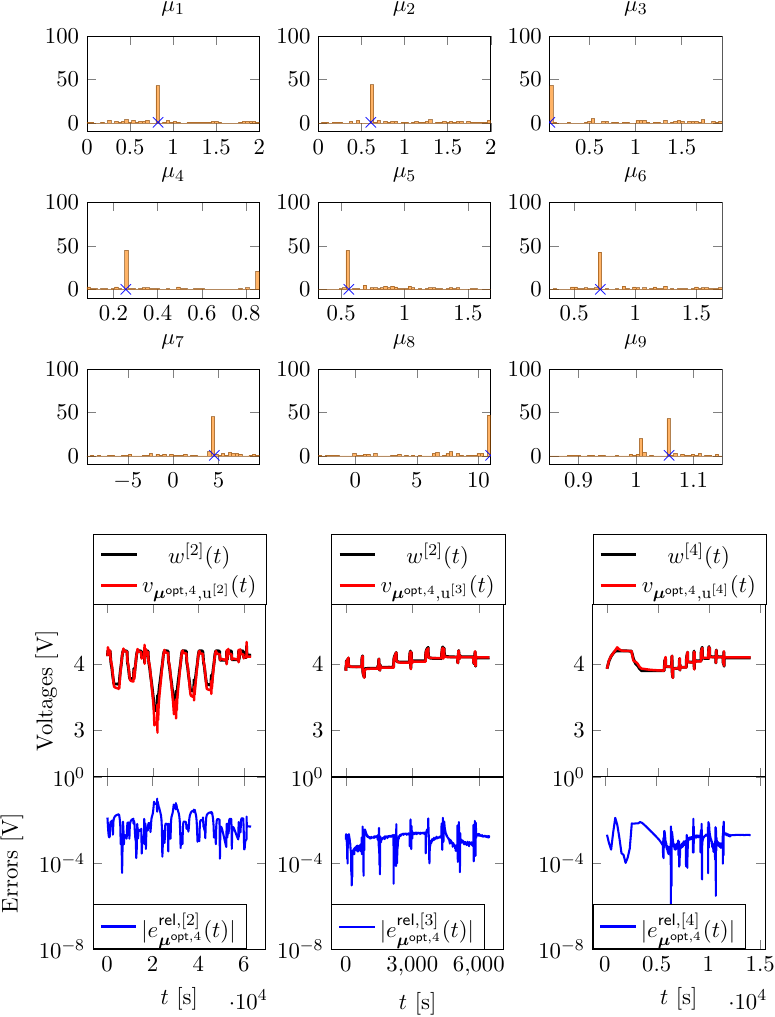}
    \caption{Results for Test 4. Top: bar plots of qualitative analysis of 100 parameter estimations starting from randomized initial guess; the blue cross finds the proposed parameter $\bmu^{\mathsf{opt},4}$; bottom: data, output and error generated by the proposed parameter.}
    \label{5:figure_test_4}
\end{figure}

As mentioned already, Test 4 is obtained by using the data of Test 3 together with the measurements available in the preparation and rest phases. The candidate parameter is definitely more concentrated around a single parameter $\bmu^{\mathsf{opt},4}$, plotted in Figure~\ref{5:figure_test_4} with the blue cross.
Counter-intuitively, the optimization over a larger set of data produces in this case a faster convergence in the parameter optimization. We interpret this effect as a result of the well-posedness of the Hessian matrix and the convexity of the parameter optimization.

We can see the results in Figure~\ref{5:figure_test_4}. The proposed parameter $\bmu^{\mathsf{opt},4}$ holds good approximation of the 3 sets of data, comparable to the one given by parameter $\bmu^{\mathsf{opt},2}$. The construction introduced in Section~\ref{Section:2.5} leads to major time cuts with the only impact being the optimization time, which averages to about \SI{154}{\second} for Test 2 and \SI{450}{\second} for Test 4, since the number of data points has almost tripled.

\subsection*{Comments on the accuracy of the optimizers}
At a first glance, it seems that we cannot find a parameter to perfectly replicate the data. This might depend on a number of factors. The first justification comes from the fact that we are working with a simplified numerical model with a reduced number of active, unknown parameters, and therefore the approximation of a real battery cell data is indeed impossible. We can also imagine that adjusting a number of parameters which now are fixed -- for example, the RK parameters -- we can fit better to the data.

On the bright side, we can make considerations about the merit of the input design. It has in fact improved robustness the parameter estimation, because it has made it harder, \ie it happens more rarely, to fall in other local optima. While it is possible that a better fitting parameter is hiding in another local optimum, it is highly unlikely.

It looks like the parameter inferred from Test 2, \ie the longest experiment, is the one better approximating the data. Test 4, which is over 4 times shorter also produces a parameter that approximates well the real data dynamics. As usual, using more expensive experiments produces data of better quality, while using shorter experiments might save time by sacrificing some of the accuracy. We believe the main advantage of using the input design algorithm is the enforcement of the well-posedness of the Hessian matrix corresponding to the parameter optimization problem. The short optimization times might be the best experimental benefit, as many optimization methods work very well with convex problems.

The input design might be even more effective when applied to more complex battery models, as it could be able to accurately represent the battery cell in all conditions: it can indeed happen that our single particle model fails at reproducing real battery cells in the settings chosen by the input design where valuable information lies, for example, for very high or very low states of charge.

We can, of course, also argue that the parameter estimation has failed because the input design was badly set up. We have selected inputs with very restricted shapes and lengths in an arbitrary and intuitive way. For example, with longer current jumps, or longer pauses, or in general longer experiment, the input design could have found the optimal parameter reproducing the data. Of course, the possibilities of generalization are endless, and we encourage future works to obtain better results.

\section{Conclusions}
In this work we have applied an algorithm for the optimal input design for the parameter estimation of a lithium-ion battery model. After introducing the numerical model, we have described a general iterative input design algorithm, where inputs are chosen in order to maximize the observability of the parameters, quantified using an approximated Fisher information matrix. By discretizing the information matrix we then define the input design problem as a finite dimensional optimization problem.

At first, we generate a collection of different input functions in order to maximize the overall observation, i.e., the collective information. The data given by this collection of input is indeed effective for the parameter estimation, and we show the effect of the input design on the convexity of the optimization problem. However, in the lab experiments, the use of the optimal inputs results in extremely long experiment times, since before and after each input some extra time needs to be taken into account, either for a resting phase or for taking the battery to a certain state of charge.

For this reason, a second input design algorithm is used, where a continuous, longer input is obtained by concatenation of optimal sub-intervals. In this way, the problem of resting and state of charge preparation is avoided -- or at least done only once -- and the overall experiment is shorter.

We derive four tests to interpret the effect of the input design and the accuracy of the parameter estimation using real battery data. In three out of four tests the input design can indeed make the parameter estimation easier, in the way that it makes us aware of a common local optimum, which we propose as the approximator of the unknown parameter.

Although we infer different parameters with different discrepancy between the numerical output and the real battery voltage, we observe how the input design was able to iteratively generate a current profile facilitating the parameter estimation, both in accuracy and in optimization time.

As mentioned, we believe that the advantage of the input design algorithms would be more evident in more complex battery cell models, where all the electro-chemical dynamics are reproduced and the data can be completely reproduced by the numerical model. In a wider scope, the techniques of optimal input/experimental design propose a goal-oriented way of defining experiment inputs and settings. Such designs might be built to guarantee several mathematical properties and rigor. The management of resources, such as time or energy, could be imposed a-priori in the design of the experiments in order to prevent waste.

\section*{Acknowledgments}
The authors would like to acknowledge the financial support within the COMET K2 Competence Centers for Excellent Technologies from the Austrian Federal Ministry for Climate Action (BMK), the Austrian Federal Ministry for Labour and Economy (BMAW), the Province of Styria (Dept.~12) and the Styrian Business Promotion Agency (SFG). The Austrian Research Promotion Agency (FFG) has been authorised for the programme management. The authors further acknowledge support by the Deutsche Forschungsgemeinschaft for the project Localized Reduced Basis Methods for PDE-constrained Parameter Optimization under contract VO 1658/6-1.

\section*{Author contributions}
\textbf{Andrea Petrocchi:} Conceptualization, Methodology, Software, Validation, Formal analysis, Writing - Original Draft, Visualization.
\textbf{Matthias K. Scharrer:} Software, Resources, Writing - Original Draft.
\textbf{Franz Pichler:} Methodology, Software.
\textbf{Stefan Volkwein:} Methodology, Writing - Review \& Editing, Supervision, Project administration, Funding acquisition.

\bibliographystyle{elsarticle-harv}

\end{document}